\documentclass[12pt,draftcls,onecolumn]{IEEEtran}

\usepackage{amssymb,amsmath,amscd}
\usepackage{graphicx}
\usepackage{eqnarray}
\usepackage{hyperref}

\newtheorem{definition}{Definition}
\newtheorem{exemple}{Example}
\newtheorem{Theorem}{Theorem}

\newtheorem{Properties}{Properties}
\newtheorem{corollaire}{Corollary}
\newtheorem{lemme}{Lemma}

\newtheorem{proposition}{Proposition}
\newtheorem{remarque}{Remark}

\newcommand{\id}{\mathsf{Id}}
\newcommand{\im}{\mathsf{Im}}
\newcommand{\z}{\mathbb{Z}}
\newcommand{\zmax}{\overline{\mathbb{Z}}_{\max}}
\newcommand{\zmaxg}{\overline{\mathbb{Z}}_{\max}[\![\gamma]\!]}

\newcommand{\ie}{ $i.e.$ }

\catcode`@=11
\def\plslash{\ifx\@currsize\normalsize
{\mathchoice
{\mbox{\raisebox{0.2ex}{$\scriptstyle\circ$}\kern-1ex$\setminus$}}
{\mbox{\raisebox{0.2ex}{$\scriptstyle\circ$}\kern-1ex$\setminus$}}%
{\mbox{\raisebox{0.14ex}{$\scriptscriptstyle\circ$}\kern-0.8ex%
${\scriptstyle\setminus}$}}%
{\,\mbox{\raisebox{0.14ex}{$\scriptscriptstyle\circ$}\kern-0.8ex%
${\scriptstyle\setminus}$}}}%
\else\ifx\@currsize\large\,\mbox{\raisebox{0.2ex}{$\scriptstyle\circ$}\kern-1ex$\setminus$}
\else\ifx\@currsize\small\,\mbox{\raisebox{0.2ex}{$\scriptstyle\circ$}\kern-1ex$\setminus$}
\else\,\mbox{\raisebox{0.2ex}{$\scriptstyle\circ$}\kern-1ex$\setminus$}
\fi\fi\fi}

\def\prslash{\ifx\@currsize\normalsize
{\mathchoice
{\mbox{\raisebox{0.2ex}{$\scriptstyle\circ$}\kern-1ex$/$}}
{\mbox{\raisebox{0.2ex}{$\scriptstyle\circ$}\kern-1ex$/$}}%
{\mbox{\raisebox{0.14ex}{$\scriptscriptstyle\circ$}\kern-0.8ex%
${\scriptstyle/}$}}%
{\mbox{\raisebox{0.14ex}{$\scriptscriptstyle\circ$}\kern-0.8ex%
${\scriptstyle/}$}}}%
\else\ifx\@currsize\large\mbox{\raisebox{0.2ex}{$\scriptstyle\circ$}\kern-1ex$/$}
\else\ifx\@currsize\small\mbox{\raisebox{0.2ex}{$\scriptstyle\circ$}\kern-1ex$/$}
\else\mbox{\raisebox{0.2ex}{$\scriptstyle\circ$}\kern-1ex$/$}
\fi\fi\fi}

\newcommand{\lslash}{\protect\plslash}
\newcommand{\rslash}{\protect\prslash}
\def\pomoins{\ifx\@currsize\normalsize
\mbox{ $\circ\kern-1.48ex-$ } \else\ifx\@currsize\large\mbox{
$\circ\kern-1.45ex-$ } \else\ifx\@currsize\small\mbox{
$\circ\kern-1.51ex-$ } \else\mbox{ $\circ\kern-1.40ex-$
}\fi\fi\fi}

\catcode`@=12

\newcommand{\rfrac}[2]{\,\rule[0.2ex]{0.1ex}{0.4ex}\kern-0.27ex%
\frac{#1\,}{\,#2}\kern-0.29ex\rule[0.55ex]{0.1ex}{0.4ex}\,\mbox{}}

\newcommand{\lfrac}[2]{\,\rule[0.55ex]{0.1ex}{0.4ex}\kern-0.29ex%
\frac{\,#1}{#2\,}\kern-0.27ex\rule[0.2ex]{0.1ex}{0.4ex}\,\mbox{}}

\begin{document}
{

\title{Observer Design for (max,plus) Linear Systems,\\
\small{Extended version of a paper published in IEEE Transaction on Automatic Control, vol.55-2, 2010, pp 538-543.}}

%
%

\author{Laurent~Hardouin\thanks{L. Hardouin, B. Cottenceau, M. Lhommeau are with the
Laboratoire d'Ing\'enierie des Syst\`emes Automatis\'es,
 Universit\'e d'Angers (France), e-mail : laurent.hardouin@univ-angers.fr,bertrand.cottenceau@univ-angers.fr,mehdi.lhommeau@univ-angers.fr},
 Carlos~Andrey~Maia\thanks{C.A. Maia is with the Departamento de Engenharia El\'{e}trica, Universidade
Federal de Minas Gerais, Av. Ant\^{o}nio Carlos 6627, Pampulha,
31270-010, Belo Horizonte, MG, Brazil, e-mail :
maia@cpdee.ufmg.br.}
 , Bertrand~Cottenceau, Mehdi~Lhommeau
}

\maketitle

\begin{abstract}

This paper deals with the state estimation for max-plus linear
systems. This estimation is carried out following the ideas of the
observer method for classical linear systems. The system matrices
are assumed to be known, and the observation of the input and of
the output is used to compute the estimated state. The observer
design is based on the residuation theory which is suitable to
deal with linear mapping inversion in idempotent semiring.
\end{abstract}


\IEEEpeerreviewmaketitle


\section{Introduction}
Many discrete event dynamic systems, such as transportation
networks \cite{farhi-05,heiderg06}, communication networks, manufacturing assembly lines \cite{Cohen-85a},
are subject to synchronization phenomena.  Timed event graphs (TEGs)
 are a subclass of timed Petri nets and are suitable tools to
model these systems. A timed event graph is a timed Petri net of
which all places have exactly one upstream transition and one
downstream transition. Its description can be transformed into a
($max,+$) or a ($min,+$) linear model and \textit{vice versa}
\cite{Cohen-84,Baccelli-92}. This property has advantaged the
emergence of a specific control theory for these systems, and
several control strategies have been proposed, $e.g.$, optimal
open loop control \cite{Cohen99,Menguy00b,lhommeau05,maia05}, and
optimal feedback control in order to solve the model matching
problem \cite{cottenceau01a,maia03b,lhommeau03a,maia05} and also
\cite{Shang05}.
 This paper focuses on observer design for ($max,+$) linear systems. The observer aims at estimating the state for a given
 plant by using input and output measurements. The state trajectories correspond to the transition
 firings of the corresponding timed event graph, their
 estimation is worthy of interest because it provides insight into internal properties of the system. For example
 these state estimations are sufficient to reconstruct the marking of the graph,
 as it is done in \cite{giua-02} for Petri nets without temporization.
  The state estimation has
 many potential applications, such as fault detection,
 diagnosis, and state feedback control.


The ($max,+$) algebra is a particular idempotent semiring,
therefore section \ref{Mplus-Syst} reviews some algebraic tools
concerning these algebraic structures. Some results about the
residuation theory and its applications over semiring are also
given. Section \ref{TEGDescription} recalls the description of
timed event graphs in a semiring of formal series. Section
\ref{Observer} presents and develops the proposed observer.
It is designed by analogy with the classical Luenberger \cite{Luenberger-71} observer
for linear systems. It is done under the assumption that the system behavior is
($max,+$)-linear. This assumption means the model represents the fastest
system behavior, in other words it implies that the
system is unable to be accelerated, and consequently the disturbances can only reduce the system performances
$i.e.$, they can only delay the events occurrence. They can be seen as machine breakdown in a manufacturing system,
or delay due to an unexpected crowd of people in a transport network.
In the opposite, the disturbances which increase system performances,
$i.e.$, which anticipate the events occurrence, could give an upper estimation of the state, in this sense the results obtained are not equivalent to the observer for the classical linear systems.
Consequently, it is assumed that the model and the initial state correspond to the fastest behavior (e.g. ideal behavior of the manufacturing system without extra delays or ideal behavior of the transport network without traffic holdup and with the maximal speed) and that disturbances only delay the occurrence of events.
Under these assumptions a sufficient condition allowing to ensure equality between the state and the estimated state is given in proposition \ref{PropositionEgalite} in spite of possible disturbances, and proposition \ref{PropositionPente} yields some weaker sufficient conditions allowing to ensure equality between the asymptotic slopes of the
state and the one of the estimated state, that means the error between  both is always bounded.
We invite the reader to consult the following link \url{http://www.istia.univ-angers.fr/~hardouin/Observer.html} to discover a dynamic illustration of the observer behavior.

\section{Algebraic Setting}
\label{Mplus-Syst}  An idempotent semiring $\mathcal{S}$ is an
algebraic structure with two internal operations denoted by
$\oplus$ and $\otimes$. The operation $\oplus$ is associative,
commutative and idempotent, that is, $a\oplus a =a$. The operation
$\otimes$ is associative (but not necessarily commutative)  and
distributive on the left and on the right with respect to
$\oplus$. The neutral elements of $\oplus$ and $\otimes$ are
represented by $\varepsilon$ and $e$ respectively, and
$\varepsilon$ is an absorbing element for the law $\otimes$
($\forall a \in \mathcal{S}, \varepsilon \otimes a = a \otimes
\varepsilon = \varepsilon ).$ As in classical algebra, the
operator $\otimes$ will be often omitted in the equations,
moreover, $a^{i} = a \otimes a^{i-1} $ and $a^{0}=e$. In this
algebraic structure, a partial order relation is defined by  $ a
\succeq b \Leftrightarrow a=a \oplus b \Leftrightarrow ~b=a \wedge b$
(where $a\wedge b$ is the greatest lower bound of $a$ and $b$),
therefore an idempotent semiring $\mathcal{S}$ is a partially
ordered set  (see \cite{Baccelli-92,heiderg06} for an exhaustive
introduction).
An idempotent semiring $\mathcal{S}$ is said to be complete if it
is closed for infinite $\oplus$-sums and if $\otimes$ distributes
over infinite $\oplus$-sums.  In particular
$\top=\bigoplus_{x\in\mathcal{S}}  x$ is the greatest
element of $\mathcal{S}$ ($\top$ is called the top element
of $\mathcal{S}$).

\begin{exemple}[$\zmax$ ]
Set $\zmax = \z \cup \lbrace -\infty,+\infty\rbrace$ endowed with
the $\max$ operator as sum and the classical sum $+$ as product is
a complete idempotent semiring, usually denoted $\zmax$, of which
$\varepsilon=-\infty$ and $e=0$.
\end{exemple}
\begin{Theorem} [see \cite{Baccelli-92}, th. 4.75]\label{Kleene}  The implicit inequality $x \succeq  ax \oplus b$
as well as the equation $x = ax \oplus b$  defined over
$\mathcal{S}$, admit $x=a^*b$ as the least solution, where $a^* =
\bigoplus\limits_{i \in \mathbb{N}} a^i$ (Kleene star operator).
\end{Theorem}

\begin{Properties}
The Kleene star operator
 satisfies the following well known properties
(see \cite{Gaubert-92} for proofs, and \cite{Krob-91} for more
general results):
\begin{eqnarray}
& a^\ast= (a^\ast)^\ast, \label{r1}& a^\ast a^\ast= a^\ast, \label{r1bis} \\
& (a \oplus b)^\ast  =  a^\ast (b a^\ast )^\ast    =  (a^\ast
b)^\ast a^\ast, &   \label{r2}  b(a b)^\ast = (ba)^\ast
b.\label{r3}
\end{eqnarray}
Thereafter, the operator $a^+=\bigoplus\limits_{i \in
\mathbb{N}^+} a^i=aa^\ast=a^\ast a$ is also considered, it
satisfies the following properties:
\begin{eqnarray}%
&a^+= (a^+)^+, \label{r3demi}  & a^\ast= e \oplus a^+, \label{r3bis}\\
& (a^\ast)^+=(a^+)^\ast = a^\ast, \label{r4}& a^+ \preceq a^\ast.
\label{r5}
\end{eqnarray}
\end{Properties}

\begin{definition}[Residual and residuated mapping]  \label{def_residuation}
An order preserving mapping $f : \mathcal{D} \rightarrow
\mathcal{E}$, where $\mathcal{D}$ and $\mathcal{E}$ are partially
ordered sets, is a \textit{residuated mapping} if for all $y \in
\mathcal{E}$ there exists a greatest solution for the inequality
$f(x) \preceq y$ (hereafter denoted $f^\sharp(y)$). Obviously, if
equality $f(x) = y$ is solvable, $f^\sharp(y)$ yields the greatest
solution. The mapping $f^\sharp$ is called the \textit{residual}
of $f$ and $f^\sharp(y)$ is the optimal solution of the
inequality.
\end{definition}

\begin{Theorem}[see \cite{Blyth-72},\cite{Baccelli-92}] \label{carac_res_ord}  Let
$f:(\mathcal{D},\preceq) \rightarrow (\mathcal{C},\preceq)$ be an
order preserving mapping. The following statements are equivalent
\begin{center}
\begin{enumerate}
\item[(i)] $f$ is residuated.
\item[(ii)] there exists an unique order preserving mapping $f^\sharp:\mathcal{C} \rightarrow \mathcal{D}$
such that $f \circ f^\sharp \preceq \id_\mathcal{C} \;\text{ and
}\; f^\sharp \circ f \succeq \id_\mathcal{D}$.
\end{enumerate}
\end{center}
\end{Theorem}

\begin{exemple}\label{exemple_La}
Mappings $\Lambda_a : x \mapsto a \otimes x$ and $\Psi_a : x
\mapsto x \otimes a$ defined over an idempotent semiring
$\mathcal{S}$ are both residuated (\cite{Baccelli-92}, p. 181).
Their residuals are order preserving mappings denoted respectively
by $\Lambda_{a}^\sharp (x) = a \lslash x $ and $\Psi_{a}^\sharp
(x) = x \rslash a $. This means that $ a \lslash b$ (resp. $b
\rslash a$) is the greatest solution of the inequality $ a \otimes
x \preceq b$ (resp. $x \otimes a \preceq b$).
\end{exemple}

\begin{definition}[Restricted mapping]
\label{definition_application_restreinte} Let $f : \mathcal{D} \to
\mathcal{C}$ be a mapping and $\mathcal{B}\subseteq \mathcal{D} $.
We will denote by $f_{|\mathcal{B}}: \mathcal{B} \to \mathcal{C}$
the mapping defined by $f_{|\mathcal{B}}=f \circ
\id_{|\mathcal{B}}$ where $\id_{|\mathcal{B}}:\mathcal{B}\to
\mathcal{D},x \mapsto x$ is the canonical injection. Identically,
let $\mathcal{E}\subseteq \mathcal{C}$ be a set such that $\im f
\subseteq \mathcal{E}$. Mapping $_{\mathcal{E}|}f : \mathcal{D}
\to \mathcal{E}$ is defined by $f=\id_{|\mathcal{E}} \circ
{_{\mathcal{E}|}f}$, where $\id_{|\mathcal{E}}: \mathcal{E} \to
\mathcal{C}, x \mapsto x$.
\end{definition}
\begin{definition}[Closure mapping]
\label{definition_fermeture} A closure mapping is an order
preserving mapping $f:\mathcal{D}\to \mathcal{D}$ defined on an
ordered set $\mathcal{D}$ such that $f\succeq \id_{\mathcal{D}}$
and $f\circ f =f$.
\end{definition}
\begin{proposition}[see \cite{cottenceau01a}]
\label{proposition_restriction_fermeture_residuable} Let
$f:\mathcal{D} \to \mathcal{D}$ be a closure mapping. Then, $_{\im
f |}f$ is a residuated mapping whose residual is the canonical
injection $\id_{|\im f}$.
\end{proposition}

\begin{exemple}\label{exempleSa}
Mapping  $K: \mathcal{S} \rightarrow \mathcal{S}, x \mapsto
x^\ast$ is a closure mapping (indeed $a \preceq a^\ast$ and
$a^\ast=(a^\ast)^\ast$ see equation (\ref{r1})). Then $(_{\im
K|}K)$ is residuated and its residual is $(_{\im K|}K)^{\sharp}=
\id_{|\im K}$. In other words, $x=a^\ast$ is the greatest solution
of inequality $x^\ast\preceq a$ if $a \in \im K$, that is $ x
\preceq a^\ast \Leftrightarrow x^\ast \preceq a^\ast$.
\end{exemple}

\begin{exemple}\label{exemplePa}
Mapping  $P: \mathcal{S} \rightarrow \mathcal{S}, x \mapsto x^+$
is a closure mapping (indeed $a \preceq a^+$ and $a^+=(a^+)^+$ see
equation (\ref{r3demi})). Then $(_{\im P|}P)$ is residuated and
its residual is $(_{\im P|}P)^{\sharp}= \id_{|\im P}$. In other
words, $x=a^+$ is the greatest solution of inequality $x^+\preceq
a$ if $a \in \im P$, that is $ x \preceq a^+ \Leftrightarrow x^+
\preceq a^+$.
\end{exemple}
\begin{remarque} \label{RemarkimKimP}
According to equation (\ref{r4}), $(a^\ast)^+ = a^\ast$, therefore
$ \im K \subset \im P$.
\end{remarque}

\begin{Properties}
Some useful results involving these residuals are presented below
(see \cite{Baccelli-92} for proofs and more complete results).
\begin{equation}
\begin{array}{ll}
 a \lslash a  =  (a \lslash a )^\ast ~~~~  & a \rslash a  =  (a \rslash a )^\ast \label{resir3}
 \end{array}
\end{equation}
\begin{equation}
\begin{array}{ll}
 a(a \lslash (a x))   = a x ~~~~  & ((xa) \rslash a) a = x a \label{resir4}
 \end{array}
\end{equation}
\begin{equation}
\begin{array}{ll}
 b \lslash a \lslash x  = (ab) \lslash x ~~~~  &
x \rslash a \rslash b  = x \rslash
(ba) \label{resir8}
 \end{array}
\end{equation}
\begin{equation}
\begin{array}{ll}
   a^\ast \lslash (a^\ast x)  =  a^\ast x ~~~~  & (a^\ast x)
\rslash a^\ast  = a^\ast x \label{resir6}
 \end{array}
\end{equation}
\begin{equation}
\begin{array}{ll}
   (a \lslash x) \wedge (a \lslash y)  =  a \lslash (x \wedge y)
~~~~ & (x \rslash a) \wedge (y \rslash a)  =  (x \wedge y) \rslash
a \label{resir9}
 \end{array}
\end{equation}

\end{Properties}

The set of $n \times n$ matrices with entries in $\mathcal{S}$ is
an idempotent semiring. The sum, the product and the residuation of
matrices are defined after the sum, the product and the residuation of
scalars in $\mathcal{S}$, \textit{i.e.},
\small
\begin{eqnarray}
&( A \otimes B  )_{ik}   =  \bigoplus\limits_{j=1 \ldots
n} ( a_{ij} \otimes b_{jk} ) &\text{~~} \label{MatrixProduct} \\
&(A \oplus B  )_{ij}  =   a_{ij} \oplus b_{ij}, &\text{~~} \label{MatrixSum} \\
&( A \lslash B )_{ij} =  \bigwedge \limits_{k=1..n} ( a_{ki} \lslash b_{kj} ) \text{~} ,\text{~}
 ( B \rslash A  )_{ij} = \bigwedge \limits_{k=1..n}  ( b_{ik} \rslash a_{jk} ).&\text{~}
\label{MatrixRightResiduation}
\end{eqnarray}
\normalsize
The identity matrix of $\mathcal{S}^{n \times n}$ is the matrix
with entries equal to $e$ on the diagonal and to $\varepsilon$
elsewhere. This identity matrix will also be denoted $e$, and the
matrix with all its entries equal to $\varepsilon$ will also be
denoted $\varepsilon$.

\begin{definition}[Reducible and irreducible
matrices]\label{definition_irreducible} Let $A$ be a $n \times n$
matrix with entries in a semiring $\mathcal{S}$. Matrix $A$ is
said reducible, if and only if for some permutation matrix $P$,
the matrix $P^T A P$ is block upper triangular. If matrix $A$ is
not reducible, it is said to be irreducible.
\end{definition}

\section{TEG description in idempotent semiring}\label{TEGDescription}
 Timed event graphs  constitute a
subclass of timed Petri nets \ie those whose places have one and
only one upstream and downstream transition. A timed event graph
(TEG) description can be transformed into a $(max,+)$ or a
$(min,+)$ linear model and \textit{vice versa}. To obtain an
algebraic model in $\zmax$, a ``dater" function is associated to
each transition. For transition labelled $x_i$,  $x_i(k)$
represents the date of the $k^{th}$ firing (see
\cite{Baccelli-92},\cite{heiderg06}). A trajectory of a TEG
transition is then a firing date sequence of this transition. This
collection of dates can be represented by a formal series $
x(\gamma)=\bigoplus_{k \in \mathbb{Z}}x_i(k)\otimes \gamma^{k}$
where $x_i(k)\in \zmax$ and $\gamma$ is a backward shift
operator\footnote{Operator $\gamma$ plays a role similar to
operator $z^{-1}$ in the $\mathcal{Z}-\text{transform}$ for the
conventional linear systems theory.} in  the event domain
(formally $\gamma x(k)= x(k-1)$). The set of formal series in
$\gamma$ is denoted by $\zmaxg$ and constitutes a complete
idempotent semiring. For instance, considering the TEG  in figure
\ref{figureGET}, daters $x_1$, $x_2$ and $x_3$ are related as
follows over $\zmax$:
$x_1(k)=4\otimes x_1(k-1) \oplus 1\otimes x_2(k) \oplus 6 \otimes
x_3(k).$ Their respective $\gamma$-transforms, expressed over
$\zmaxg$, are then related as:
$$x_1(\gamma)=4 \gamma x_1(\gamma) \oplus 1 x_2(\gamma) \oplus
6 x_3(\gamma).$$ In this paper  TEGs are modelled in this
setting, by the following  model  :
\begin{eqnarray}\label{RealModel}
  \nonumber x & = & Ax \oplus Bu \oplus Rw\\
  y & = & Cx,
\end{eqnarray}

\noindent where $u \in (\zmaxg)^p $, $y \in (\zmaxg)^m$ and $x \in
(\zmaxg)^n $ are respectively the controllable input, output and
state vector, $i.e.$, each of their entries is a trajectory which
represents the collection of firing dates of the corresponding
transition.
 Matrices $A \in (\zmaxg)^{n \times n}$, $B \in (\zmaxg)^{n \times p}$, $C
\in (\zmaxg)^{m \times n}$ represent the links between each
transition, and then describe the structure of the graph. Vector
$w \in (\zmaxg)^l$ represents uncontrollable inputs ($i.e.$
disturbances\footnote{In manufacturing setting, $w$ may represent
machine breakdowns or failures in component supply.}). Each entry
of $w$ corresponds to a transition which disables the firing of
internal transition of the graph, and then decreases the
performance of the system. This vector is bound to the graph
through matrix $R \in (\zmaxg)^{n \times l}$.

Afterwards, each input transition  $u_i$ (respectively $w_i$) is
assumed to be connected to one and only one internal transition
$x_j$, this means  that each column of matrix $B$ (resp. $R$) has
one entry equal to $e$ and the others equal to $\varepsilon$ and
at  most one entry equal to $e$ on each row. Furthermore, each
output transition $y_i$ is assumed to be linked to one and only
one internal transition $x_j$, $i.e$ each row of matrix $C$ has
one entry equal to $e$ and the others equal to $\varepsilon$ and
at  most one entry equal to $e$ on each column. These requirements
are satisfied without loss of generality, since it is sufficient
to add extra input and output transition. Note that if $R$ is
equal to the identity matrix, $w$ can represent initial state of
the system $x(0)$ by considering $w=x(0)\gamma^0 \oplus...$ (see
\cite{Baccelli-92}, p. 245, for a discussion about compatible
initial conditions).
By considering theorem \ref{Kleene}, this system can be rewritten
as :
\begin{eqnarray}
   x & = &A^\ast B u \oplus A^\ast R w \label{RealStatetransfer1}\\
   y& = &CA^\ast B u \oplus CA^\ast R w , \label{RealStatetransfer2}
  \end{eqnarray}

where $(CA^\ast B) \in (\zmax)^{m\times p}$ (respectively $
(CA^\ast R) \in (\zmax)^{m\times l}$) is the input/output (resp.
disturbance/output) transfer matrix.
Matrix $(CA^\ast B)$ represents the earliest behavior of the
system, therefore it must be underlined that the uncontrollable
inputs vector $w$ (initial conditions or disturbances) is only
able to delay the transition firings, \ie, according to the order
relation of the semiring, to increase the vectors $x$ and $y$.

If the TEG is strongly connected, \ie there exists at least one path
between transitions $x_i, x_j~\forall i,j$, then  matrix $A$ is
irreducible. If $A$ is reducible, according to definition
\ref{definition_irreducible}, there exists a permutation matrix such
that :
\begin{eqnarray}
A=\begin{pmatrix} A_{11} & A_{12}&...& A_{1k} \\
\varepsilon & A_{22}&...& A_{2k} \\
\vdots & \vdots &\ddots& \vdots \\
\varepsilon & \varepsilon & ...& A_{kk}
\end{pmatrix}
\end{eqnarray}

where $k$ is the number of strongly connected components of the
TEG, and each matrix $A_{ii}$ is an irreducible matrix associated
to the component $i$. Matrices $A_{ij}$ (with $i \neq j$)
represent the links between these strongly connected components.
\begin{figure}[h]
\begin{center}
\includegraphics[scale=0.45]{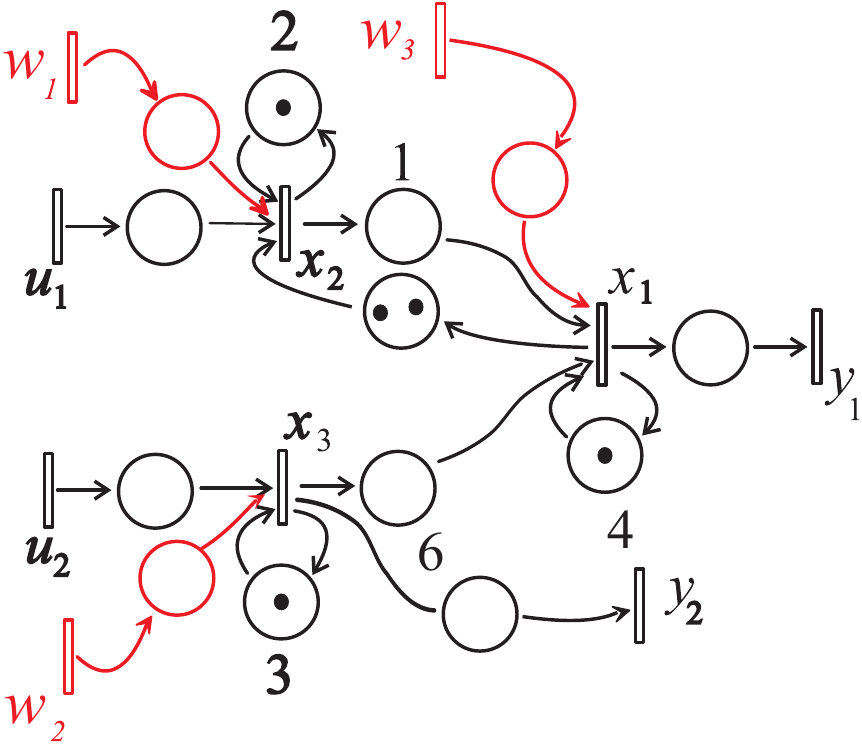}
\end{center}
\caption{\label{figureGET}Timed event graph, $u_i$ controllable  and  $w_i$ uncontrollable inputs.}
\end{figure}
 Consequently, for the TEG depicted fig.\ref{figureGET}, the following matrices are obtained:
 \small
$A=\begin{pmatrix} 4\gamma & 1 & 6\\
 \gamma^2 & 2\gamma &  \varepsilon\\
 \varepsilon & \varepsilon & 3\gamma
\end{pmatrix},\; B= \begin{pmatrix} \varepsilon & \varepsilon\\ e & \varepsilon \\ \varepsilon &
e
\end{pmatrix},\;$ $ C=\begin{pmatrix}e & \varepsilon & \varepsilon \\
  \varepsilon & \varepsilon & e
\end{pmatrix}, \;R=\begin{pmatrix}e & \varepsilon & \varepsilon \\\varepsilon & e &
\varepsilon\\ \varepsilon & \varepsilon & e
\end{pmatrix},$
\normalsize
and leads to the following $A^\ast$ matrix:
\small
\begin{eqnarray}\label{Astar}
A^\ast  = \begin{pmatrix} (4\gamma)^\ast & 1(4\gamma)^\ast
& 6(4\gamma)^\ast\\
\gamma^2(4\gamma)^\ast & e \oplus 2\gamma \oplus 4 \gamma^2\oplus 6 \gamma^3 \oplus 9\gamma^4(4\gamma)^\ast & 6\gamma^2(4\gamma)^\ast \\
\varepsilon &  \varepsilon & (3\gamma)^\ast \end{pmatrix} \nonumber
\end{eqnarray}
\normalsize
According to assumptions about matrices $C$, $B$, and $R$, the
 matrices $(CA^\ast B)$ and $(CA^\ast R)$
are composed of some  entries of matrix $A^\ast$. Each entry
is a periodic series \cite{Baccelli-92} in the $\zmaxg$ semiring.
A periodic series $s$ is usually represented by $s=p\oplus q
r^\ast$, where $p$ (respectively $q$) is a polynomial depicting
the transient (resp. the periodic) behavior, and
$r=\tau\gamma^\nu$ is a monomial depicting the periodicity
allowing to define the asymptotic slope of the series as
$\sigma_{\infty}(s)=\nu/\tau$ (see figure \ref{periodicseries}).
Sum, product, and residuation of periodic series are well defined
(see \cite{Gaubert-92}), and algorithms and software toolboxes are
available in order to handle periodic series and compute transfer
relations (see \cite{lhommeau00a}). Below, only the rules between
monomials and properties concerning asymptotic slope are recalled
:
\begin{equation}
\begin{array}{lcl} \nonumber
t_1\gamma^n \oplus t_2 \gamma^n&=&max(t_1,t_2)\gamma^n,\\
t\gamma^{n_1} \oplus t \gamma^{n_2}&=&t\gamma^{min(n_1,n_2)},\\
t_1\gamma^{n_1} \otimes t_2 \gamma^{n_2}&=&(t_1+t_2)\gamma^{(n_1+n_2)},\\
(t_1\gamma^{n_1}) \rslash (t_2 \gamma^{n_2})&=&(t_2\gamma^{n_2}) \lslash (t_1 \gamma^{n_1})=(t_1-t_2)\gamma^{(n_1-n_2)},
\end{array}
\end{equation}
\begin{eqnarray}
\sigma_{\infty}(s \oplus s') &=&
\min(\sigma_{\infty}(s),\sigma_{\infty}(s')), \label{penteoplus}\\
\sigma_{\infty}(s \otimes  s') &=&
\min(\sigma_{\infty}(s),\sigma_{\infty}(s')),
\label{penteotimes}\\
 \sigma_{\infty}(s \wedge  s') &=&
\max(\sigma_{\infty}(s),\sigma_{\infty}(s')), \label{pentewedge}
\end{eqnarray}
\begin{equation}
\begin{array}{lll}
\textnormal{if~} \sigma_{\infty}(s) \leq \sigma_{\infty}(s')
\textnormal{~then~} \sigma_{\infty}(s' \lslash s)& = &
\sigma_{\infty}(s),\\ \textnormal{~else~~} s' \lslash s=\varepsilon.
\label{penteresiduation}
\end{array}
\end{equation}
\begin{figure}[h]
\begin{center}
\includegraphics[scale=0.55]{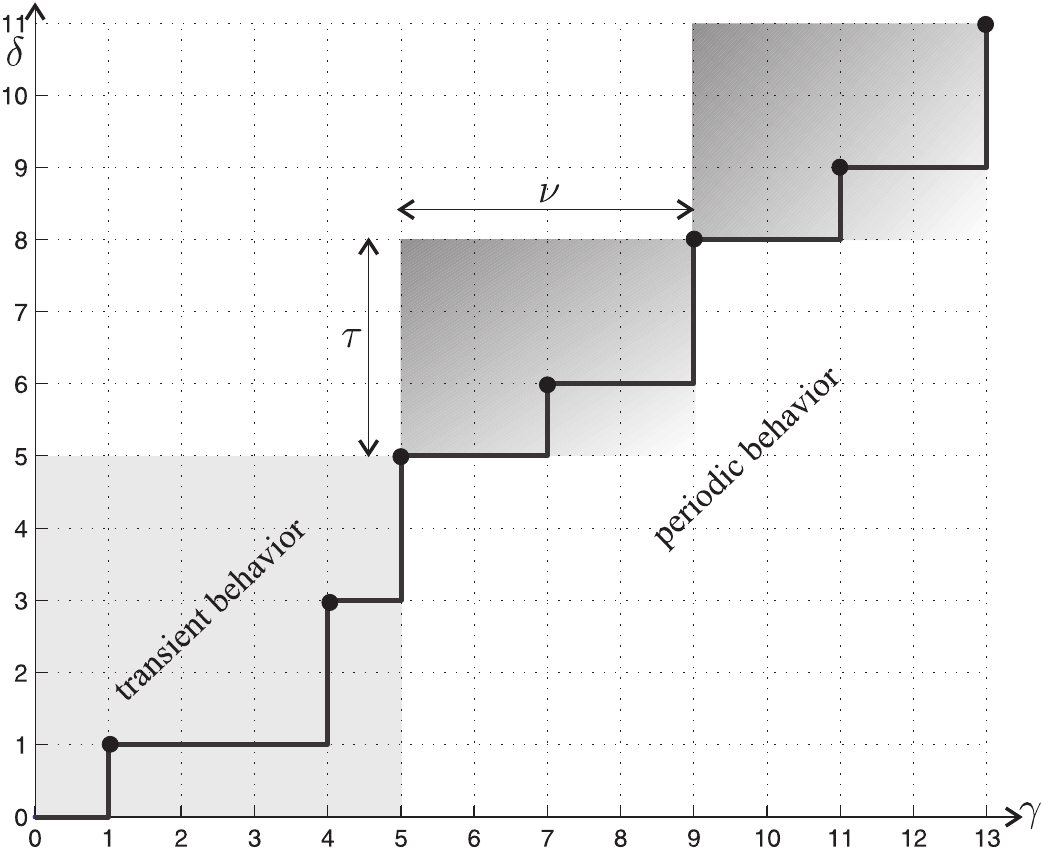}
\end{center}
\caption{\label{periodicseries}Periodic series $s=(e \oplus 1
\gamma^1 \oplus 3 \gamma^4) \oplus (5\gamma^5 \oplus 6
\gamma^7)(3\gamma^4)^\ast$.}
\end{figure}

Let us recall that if matrix $A$ is irreducible then all the entries of
matrix $A^\ast$ have the same asymptotic slope, which will be
denoted $\sigma_{\infty}(A)$. If $A$ is a reducible matrix assumed
to be in its  block upper triangular representation, then
matrix $A^\ast$ is  block upper triangular and matrices
$(A^\ast)_{ii}$ are such that $(A^\ast)_{ii}=A^\ast_{ii}$ for each
$i \in [1,k]$. Therefore, since $A_{ii}$ is irreducible, all the
entries of matrix $(A^\ast)_{ii}$ have the same asymptotic slope
$\sigma_{\infty}((A^\ast)_{ii})$. Furthermore, entries of each
matrix $(A^\ast)_{ij}$ with $i<j$ are such that their asymptotic
slope is lower than or equal to
$min(\sigma_{\infty}((A^\ast)_{ii}),\sigma_{\infty}((A^\ast)_{jj}))$.

\section{Max-plus Observer}  \label{Observer}
\begin{figure}[h]
\begin{center}
\includegraphics[scale=0.6]{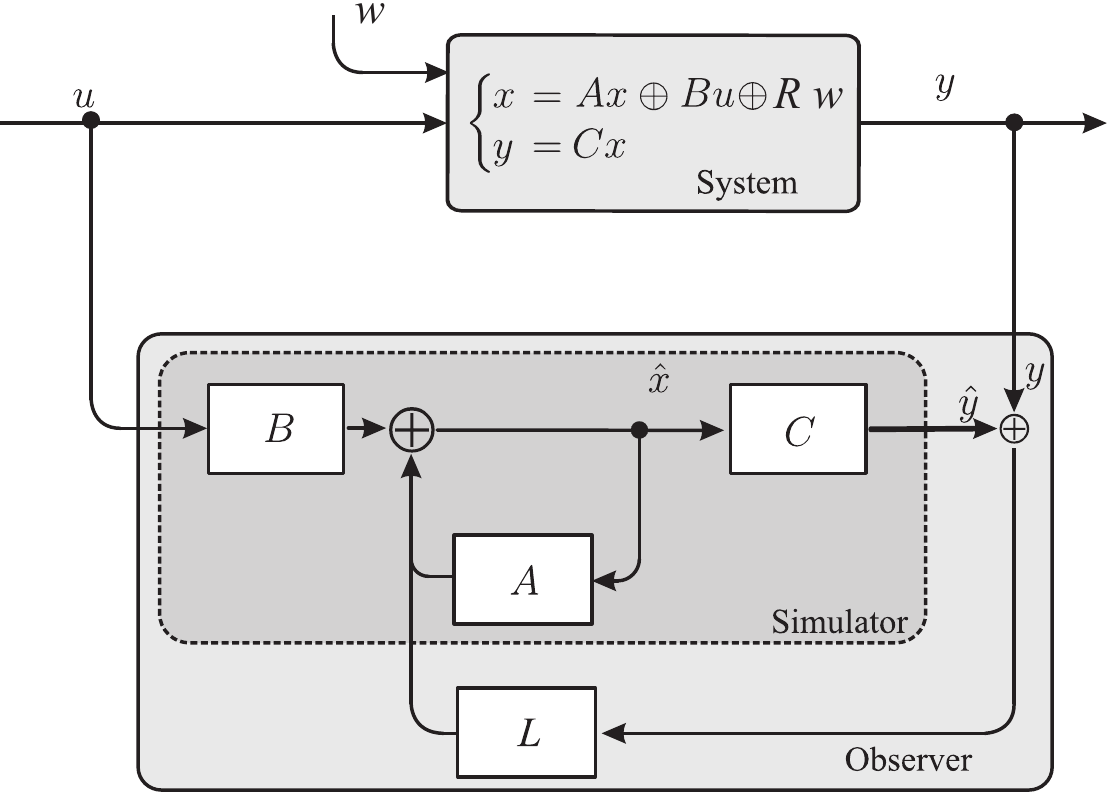}
\end{center}
\caption{\label{Obs-Syst}Observer structure.}
\end{figure}
Figure \ref{Obs-Syst} depicts the observer structure directly
inspired from the classical linear system theory (see
\cite{Luenberger-71}). The observer matrix $L$ aims at providing
information from the system output into the simulator, in order to
take  the disturbances $w$ acting on the system into account. The
simulator is described by the model\footnote{Disturbances are
uncontrollable and \textit{a priori} unknown, then the simulator
does not take them into account.} (matrices $A$, $B$, $C$) which
is assumed to represent the fastest behavior of the real system in
a guaranteed way\footnote{Unlike in the conventional linear system
theory, this assumption means that the fastest behavior of the
system is assumed to be known and that the disturbances can only
delay its behavior.}, furthermore the simulator is initialized by
the canonical initial conditions (\ie $\hat{x_i}(k)=\varepsilon,
\forall k\leq 0$). These assumptions induce that $y \succeq
\hat{y}$ since disturbances and initial conditions, depicted by
$w$, are only able to increase the system output. By considering
the configuration of figure \ref{Obs-Syst} and these assumptions,
the computation of the optimal observer matrix $L_x$ will be
proposed in order to achieve  the constraint $\hat{x} \preceq x$.
Optimality means that the matrix is obtained thanks to the
residuation theory and then it is the greatest one (see definition
\ref{def_residuation}), hence the estimated state $\hat{x}$ is the
greatest which achieves the objective. Obviously this optimality
is only ensured under the assumptions considered ($i.e.$
$\hat{y}\preceq y$).
 As in the development proposed in conventional linear systems theory, matrices $A$, $B$, $C$ and $R$
are assumed to be known, then the system transfer is given by
equations (\ref{RealStatetransfer1}) and
(\ref{RealStatetransfer2}). According to figure \ref{Obs-Syst} the
observer equations are given by:
\begin{eqnarray}\label{ObserverModel}
\nonumber   \hat{x}&=& A \hat{x}\oplus Bu \oplus L(\hat{y}\oplus y)\\
                    &=& A \hat{x}\oplus Bu \oplus LC\hat{x}\oplus LC x\\
 \nonumber  \hat{y}&=&C\hat{x}.
\end{eqnarray}
By applying Theorem
\ref{Kleene} and by considering equation
(\ref{RealStatetransfer1}), equation (\ref{ObserverModel}) becomes
:
\begin{eqnarray}
  \hat{x}&=& (A \oplus LC)^\ast Bu  \oplus (A \oplus
LC)^\ast LCA^\ast B u \nonumber \\
 & & \oplus (A \oplus LC)^\ast LCA^\ast Rw. \label{ObserverModelbis}
  \end{eqnarray}
By applying equation (\ref{r2}) the following equality is obtained :
\begin{eqnarray}\label{equalityALC}
 (A \oplus LC)^\ast=A^\ast (LC A^\ast)^\ast,
  \end{eqnarray}
by replacing in equation (\ref{ObserverModelbis}) :
\begin{eqnarray}
  \nonumber \hat{x}&=& A^\ast(LC A^\ast)^\ast Bu \oplus A^\ast (LCA^\ast )^\ast LCA^\ast B u  \\
 \nonumber & & \oplus A^\ast (LCA^\ast )^\ast LCA^\ast
   Rw,
       \end{eqnarray}
and by recalling that $(LCA^\ast)^\ast LCA^\ast=(LCA^\ast)^+$,
this equation may be written as follows :
\begin{eqnarray}
\nonumber   \hat{x} &=& A^\ast(LC A^\ast)^\ast Bu \oplus A^\ast
(LCA^\ast )^+ B u  \oplus A^\ast (LCA^\ast )^+ Rw.
\end{eqnarray}
Equation (\ref{r5}) yields $(LCA^\ast)^\ast \succeq (LCA^\ast)^+$,
then the observer model may be written as follows :
\begin{eqnarray}\label{ObserverModelSimplifie}
\nonumber    \hat{x}&=& A^\ast(LC A^\ast)^\ast Bu \oplus  A^\ast (LCA^\ast )^+ Rw\\
   &=& (A \oplus LC)^\ast Bu  \oplus (A \oplus LC)^\ast LCA^\ast Rw.
\end{eqnarray}
As said previously the objective considered is to compute the
greatest observation matrix $L$ such that the estimated state
vector $\hat{x}$ be as close as possible to state $x$, under the
constraint $\hat{x}\preceq x$, formally it can be written :
$$
\begin{array}{lcll}
(A \oplus LC)^\ast Bu \oplus  (A \oplus LC)^\ast LCA^\ast Rw & \preceq & A^\ast Bu \oplus A^\ast  Rw \hspace{1cm} & \forall (u,w) \\
\end{array}
    $$
or equivalently :
\begin{eqnarray}\label{ObjectiveObserver}
 (A \oplus L C)^\ast B  &\preceq&  A^\ast B   \\
  (A \oplus LC)^\ast LCA^\ast R &\preceq& A^\ast R .
    \end{eqnarray}
\begin{lemme} \label{LemmaL1} The greatest matrix $L$ such that $(A \oplus LC)^*B =
 A^*B$ is given by:
\begin{equation} \label{Eq-observ1}
L_ {1}= ( A^\ast B) \rslash (CA^\ast B).
\end{equation}
\end{lemme}
\begin{IEEEproof}
First let us note that $L= \varepsilon \in \zmax^{n \times m} $ is
a solution, indeed  $(A \oplus \varepsilon C)^ \ast B = A^\ast B$.
Consequently, the greatest solution of the inequality $(A \oplus
LC)^\ast B \preceq A^\ast B$ will satisfy the equality.
Furthermore, according to equation (\ref{r3}), $(A \oplus LC)^*B=
(A^\ast LC)^\ast A^\ast B$. So the objective is given by :
$$
\begin{array}{lll}
\nonumber &(A^*LC)^* A^*B \preceq A^*B \\
\nonumber \Leftrightarrow &(A^*LC)^* \preceq (A^*B) \rslash (A^*B)
& \textnormal{(see  example \ref{exemple_La})} \\
\nonumber \Leftrightarrow & (A^*LC)^* \preceq((A^*B)\rslash (A^*B))^* & \textnormal{(see eq.(\ref{resir3}))} \\
\nonumber  \Leftrightarrow & (A^*LC) \preceq (A^*B) \rslash (A^*B)
&  \textnormal{(see
example \ref{exempleSa})}\\
\nonumber \Leftrightarrow & L
 \preceq A^\ast \lslash (A^\ast B)
\rslash
(A^\ast B) \rslash C  & \textnormal{(see  example \ref{exemple_La})} \\
\nonumber \Leftrightarrow & L \preceq A^\ast \lslash (A^\ast B)\rslash (CA^\ast B) & \textnormal{(see eq.(\ref{resir8}}))\\
\Leftrightarrow & L \preceq  (A^\ast B) \rslash (CA^\ast B) =
L_{1} & \textnormal{(see eq.(\ref{resir6}))}
\end{array}
$$
\end{IEEEproof}
\begin{lemme} \label{LemmaL2} The greatest matrix $L$ that satisfies $(A \oplus LC)^\ast LCA^\ast  R \preceq A^\ast R$ is given by:
\begin{equation} \label{Eq-observ2}
L_ {2}=  ( A^\ast R) \rslash (C A^\ast R).
\end{equation}
\end{lemme}
\begin{IEEEproof}
$$
\begin{array}{ll}
\nonumber &(A \oplus LC)^\ast LCA^\ast R  \preceq A^\ast R \\
\nonumber \Leftrightarrow &A^\ast (LCA^\ast)^\ast  LCA^\ast R
\preceq A^\ast R  \textnormal{~~~~~~~~~~~~~~~~(see eq.(\ref{equalityALC}))}, \\
\nonumber \Leftrightarrow & (LCA^\ast)^\ast LCA^\ast R
\preceq A^\ast \lslash (A^\ast R)= A^\ast R\\
\nonumber & \textnormal{(see example \ref{exemple_La} and eq.(\ref{resir6}), with } x=R),\\
\nonumber \Leftrightarrow & (LCA^\ast)^\ast LCA^\ast A^\ast
R=(LCA^\ast)^+ A^\ast R \preceq
A^\ast R\\
\nonumber &  \textnormal{(see eq.(\ref{r1bis}) and $a^+$ definition)},\\
 \nonumber \Leftrightarrow &
(LCA^\ast)^+\preceq (A^\ast R) \rslash (A^\ast R)=((A^\ast R) \rslash (A^\ast R))^\ast\\
\nonumber &  \textnormal{(see  eq.(\ref{resir3}))},\\
\end{array}
$$
according to remark \ref{RemarkimKimP} the right member is in $\im
P$, then by applying the result presented in example
\ref{exemplePa}, this inequality may be written as follows :
$$
\begin{array}{lllcl}
& LCA^\ast\preceq (A^\ast R) \rslash (A^\ast R) \\
\Leftrightarrow & L\preceq (A^\ast R) \rslash (A^\ast R) \rslash
(CA^\ast )= (A^\ast R) \rslash (C A^\ast A^\ast R)\\
& \textnormal{(see example \ref{exemple_La} and eq. (\ref{resir6}))}\\
\Leftrightarrow & L\preceq (A^\ast R) \rslash (C  A^\ast R)=L_2\\
&\textnormal{(see eq. (\ref{r1bis}))}.
\end{array}
$$

\end{IEEEproof}
\begin{proposition} \label{PropL} $L_x= L_1 \wedge L_2$ is the greatest observer matrix such that:
$$\hat{x}=A \hat{x}\oplus Bu \oplus L(\hat{y}\oplus y) \preceq x=A
x \oplus Bu \oplus Rw ~~~ \forall (u,w).$$
\end{proposition}
\begin{IEEEproof}
 Lemma \ref{LemmaL1} implies $L\preceq L_1$   and lemma $\ref{LemmaL2}$ implies $L\preceq
 L_2$, then $L\preceq L_1 \wedge L_2=L_x$.
\end{IEEEproof}
\begin{corollaire}\label{CorollaryOutputEgality}
The matrix $L_x$ ensures the equality between estimated output
$\hat{y}$ and measured output $y$, $i.e.$
\begin{eqnarray}
C(A \oplus L_xC)^\ast B   &= &     CA^\ast B, \label{egalite1Lx}\\
C(A \oplus L_xC)^\ast L_xCA^\ast R  &= &CA^\ast
R.\label{egalite2Lx}
\end{eqnarray}
\end{corollaire}
\begin{IEEEproof}
Let $\tilde{L}=e \rslash C$ be a particular observer matrix.
Definition \ref{def_residuation} yields $\tilde{L}C\preceq e$ then
$(A\oplus \tilde{L}C)^\ast = A^\ast$. This equality implies
 $(A \oplus \tilde{L}C)^\ast B  =  A^\ast B$,
therefore according to lemma \ref{LemmaL1} $\tilde{L}\preceq L_1$,
since $L_1$ is the greatest solution.  That implies also that
$L_1$ is solution of equation (\ref{egalite1Lx}).
 Equality $(A\oplus \tilde{L}C)^\ast = A^\ast$ and inequality $\tilde{L}C\preceq e$ yield  $(A \oplus
\tilde{L}C)^\ast \tilde{L}CA^\ast R  =A^\ast \tilde{L}CA^\ast
R\preceq A^\ast R$ then according to lemma \ref{LemmaL2}
$\tilde{L}\preceq L_2$ since $L_2$ is the greatest solution. That
implies also that $\tilde{L}$ and $L_2$ are such that $C(A \oplus
\tilde{L}C)^\ast \tilde{L}CA^\ast R  \preceq C(A \oplus L_2C)^\ast
L_2CA^\ast R
\preceq CA^\ast R $.
 The assumption about matrix $C$ (see section \ref{TEGDescription})
yields $CC^T=e$ and $\tilde{L}=e\rslash C=C^T$, therefore $C(A
\oplus \tilde{L}C)^\ast \tilde{L}CA^\ast R=CA^\ast \tilde{L}
CA^\ast R=(C \tilde{L} \oplus CA\tilde{L} \oplus ...)CA^\ast R
\succeq C \tilde{L} CA^\ast R=C C^T CA^\ast R =CA^\ast R$.
Therefore, since  $\tilde{L}\preceq L_2$, we have $C(A \oplus
\tilde{L}C)^\ast \tilde{L}CA^\ast R=C(A \oplus L_2 C)^\ast
L_2CA^\ast R=CA^\ast R$ and both $\tilde{L}$ and $L_2$ yield
equality (\ref{egalite2Lx}). To conclude  $\tilde{L} \preceq
L_{1} \wedge L_{2}=L_x$, hence, $L_x\preceq L_1$ yields the
equality (\ref{egalite1Lx}) and $L_x\preceq L_2$ yields
(\ref{egalite2Lx}). Therefore equality $\hat{y} =y$ is ensured.
\end{IEEEproof}
\begin{remarque}\label{remarqueBbarre}
By considering matrix $\overline{B}=\begin{pmatrix}B & R
\end{pmatrix}$, equations (\ref{MatrixRightResiduation}) and
(\ref{resir9}), matrix $L_x$ may be written as : $ L_x =(A^\ast
\overline{B} \lslash (CA^\ast \overline{B})$.

According to the residuation theory (see definition
\ref{def_residuation}), $L_x$ yields $x=\hat{x}$ if possible.
Nevertheless, two questions arise, firstly is it possible to ensure
equality between the asymptotic slope of each state vector entries ?
Secondly is it possible to ensure equality between these vectors ?
Below,  sufficient conditions allowing to answer positively are
given.
\end{remarque}
\begin{proposition}\label{PropositionPente}
Let $k$ be the number of strongly connected components of the TEG
considered.
 If matrix $C \in \zmaxg^{k \times n}$ is defined as in section \ref{TEGDescription} and such that each
strongly connected component is linked to one and only one output
then $\sigma_{\infty}(x_i)=\sigma_{\infty}(\hat{x}_i) \forall i
\in [1,n] $.
\end{proposition}
\begin{IEEEproof}
 First, assuming that matrix $A$ is irreducible ($\ie$, $k=1$), then all entries
of matrix $A^\ast$ have the same asymptotic slope
$\sigma_{\infty}(A^\ast)$. As said in section \ref{TEGDescription}
entries of matrices $B$, $R$, and $C$ are equal to $\varepsilon$
or $e$, therefore, according to matrices operation definitions
(see equations (\ref{MatrixProduct}) to
(\ref{MatrixRightResiduation}) and rules (\ref{penteoplus}) to
(\ref{penteresiduation})), all the entries of matrices $A^\ast B$,
$A^\ast R$, $CA^\ast B$, $CA^\ast R$ and $L_x$ have the same
asymptotic slope which is equal to $\sigma_{\infty}(A^\ast)$.
Consequently,  by considering equation
(\ref{ObserverModelSimplifie}), $\sigma_{\infty}(((A \oplus
L_xC)^\ast B)_{ij})=\sigma_{\infty}((A^\ast B)_{ij})$ and
$\sigma_{\infty}(((A \oplus L_xC)^\ast L_xCA^\ast
R)_{ij})=\sigma_{\infty}((A^\ast R)_{ij})$, which leads to
$\sigma_{\infty}(x_i)=\sigma_{\infty}(\hat{x}_i)
\textnormal{~}\forall i \in [1,n]$.\\
Now the reducible case is considered. To increase the readability,
matrices $B$ and $R$ are assumed to be equal to $e$ and the proof
is given for a graph with two strongly connected components. The
extension for a higher dimension may be obtained in an analogous
way. As said in section \ref{TEGDescription}, matrix $A^\ast$
is block upper diagonal :
\begin{eqnarray} \nonumber
A^\ast=\begin{pmatrix} (A^\ast)_{11} & (A^\ast)_{12}\\
\varepsilon & (A^\ast)_{22}\\
\end{pmatrix},
\end{eqnarray}
all the entries of the square matrix $(A^\ast)_{ii}$ have the same
asymptotic slope $\sigma_\infty((A^\ast)_{ii})$ and all the
entries of matrix $(A^\ast)_{12}$ have the same asymptotic slope,
$\sigma_{\infty}((A^\ast)_{12})=min(\sigma_{\infty}((A^\ast)_{11}),\sigma_{\infty}((A^\ast)_{22}))$.
Assumption about matrix $C \in \zmaxg^{2 \times n}$, $i.e.$ one and
only one entry is linked to each strongly connected component,
yields the following block upper diagonal matrix :

\begin{eqnarray} \nonumber
CA^\ast=\begin{pmatrix} (CA^\ast)_{11} & (CA^\ast)_{12}\\
\varepsilon & (CA^\ast)_{22}\\
\end{pmatrix},
\end{eqnarray}
where $\begin{pmatrix} (CA^\ast)_{11} &
(CA^\ast)_{12}\end{pmatrix} $ is one row of matrix
$\begin{pmatrix} (A^\ast)_{11} & (A^\ast)_{12}\end{pmatrix}$ and
$\begin{pmatrix} \varepsilon & (CA^\ast)_{22}\end{pmatrix} $ is
one row of  matrix $\begin{pmatrix} \varepsilon &
(A^\ast)_{22}\end{pmatrix}$, hence
$\sigma_\infty((CA^\ast)_{ij})=\sigma_\infty((A^\ast)_{ij})$. Matrix $L_x$ is also block upper diagonal :

\begin{eqnarray} \nonumber
L_x=A^\ast \rslash C=\begin{pmatrix} L_{x11} & L_{x12}\\
\varepsilon & L_{x22}\\
\end{pmatrix},
\end{eqnarray}

where $\begin{pmatrix} L_{x11} & \varepsilon\end{pmatrix}^T $ is
one column of  matrix $\begin{pmatrix} (A^\ast)_{11} &
\varepsilon\end{pmatrix}^T$ and $\begin{pmatrix} L_{x12} &
L_{x22}\end{pmatrix}^T $ is one column of  matrix
$\begin{pmatrix} (A^\ast)_{12} & (A^\ast)_{22}\end{pmatrix}^T$,
hence $\sigma_\infty(L_{xij})=\sigma_\infty((A^\ast)_{ij})$.
Therefore $L_xCA^\ast$ is block upper diagonal  :
\small
\begin{eqnarray} \nonumber
L_xCA^\ast & =&\begin{pmatrix} L_{x11}(CA^\ast)_{11} & L_{x11}(CA^\ast)_{12} \oplus L_{x12}(CA^\ast)_{22}\\
\varepsilon & L_{x22}(CA^\ast)_{22}
\end{pmatrix}\\& = &  \begin{pmatrix} (L_xCA^\ast)_{11} & (L_xCA^\ast)_{12}\\
\varepsilon & (L_xCA^\ast)_{22}
\end{pmatrix},
\end{eqnarray}
\normalsize
and by considering rules (\ref{penteoplus}) and
(\ref{penteresiduation}), the sub matrices are such that
$\sigma_\infty((L_xCA^\ast)_{ij})=\sigma_\infty((A^\ast)_{ij})$.
By recalling that $(A\oplus L_xC)^\ast= A^\ast (L_xCA^\ast)^\ast$,
we obtain $\sigma_{\infty}(((A \oplus L_xC)^\ast
)_{ij})=\sigma_{\infty}((A^\ast )_{ij})$ and $\sigma_{\infty}(((A
\oplus L_xC)^\ast L_xCA^\ast )_{ij})=\sigma_{\infty}((A^\ast
)_{ij})$, which leads to
$\sigma_{\infty}(x_i)=\sigma_{\infty}(\hat{x}_i)
\textnormal{~}\forall i \in [1,n]$.
\end{IEEEproof}
\begin{proposition} \label{PropositionEgalite}
If matrix  $A^\ast \overline{B}$ is in $\im \Psi_{CA^\ast
\overline{B}}$,  matrix $L_x$ is such that $\hat{x}=x$.
\end{proposition}
\begin{IEEEproof}
First, let us recall that $$
\begin{array}{lll}
A^\ast \overline{B} \in \im
\Psi_{CA^\ast \overline{B}} & \Leftrightarrow & \exists z \textnormal{
s.t. } A^\ast \overline{B}= z CA^\ast \overline{B} \\ & \Leftrightarrow &
((A^\ast \overline{B})\rslash (CA^\ast \overline{B})) (CA^\ast
\overline{B})=  A^\ast \overline{B}.
\end{array}
$$
If  $\exists z \textnormal{
s.t. } A^\ast \overline{B}= z CA^\ast \overline{B}$ then
$$\begin{array}{lll}
 L_xCA^\ast \overline{B} &= &((A^\ast \overline{B}) \rslash (CA^\ast \overline{B}))CA^\ast \overline{B} \\
 & = &((zCA^\ast \overline{B}) \rslash (CA^\ast \overline{B})) CA^\ast \overline{B} \\
 & = &zCA^\ast \overline{B}=A^\ast \overline{B} \textnormal{  (see eq.
(\ref{resir4}),}
\end{array}
$$
by recalling that $\overline{B}=\begin{pmatrix}B &
R\end{pmatrix}$, this equality can be written
$$\begin{pmatrix}  L_xCA^\ast B & L_xCA^\ast R \end{pmatrix}= \begin{pmatrix}  A^\ast B & A^\ast R
\end{pmatrix}.
$$
Therefore $(A \oplus L_x C)^\ast L_x CA^\ast R =A^\ast (L_x C
A^\ast )^\ast L_x CA^\ast R = A^\ast (L_x CA^\ast)^+ R  = A^\ast
(L_xCA^\ast R \oplus
 (L_xCA^\ast)^2 R \oplus (L_xCA^\ast)^3 R \oplus ...)$ (see equation (\ref{equalityALC}) and $a^+$ definition).
   Since $ L_xCA^\ast R = A^\ast R$, the following equality is satified  $(L_xCA^\ast)^2 R=L_xCA^\ast A^\ast
R=L_xCA^\ast R=A^\ast R$ and more generally $(L_xCA^\ast )^i R=
A^\ast R$, therefore $L_x$ ensures equality $(A \oplus L_x C)^\ast
L_x CA^\ast R=A^\ast (L_x CA^\ast)^+ R=A^\ast R$. On the other
hand lemma \ref{LemmaL1} yields the equality $(A \oplus L_xC)^*B =
 A^*B$, which concludes the proof.
\end{IEEEproof}
\begin{remarque}
This sufficient condition gives an interesting test to know if the
number of sensors  is sufficient and if they are well localized to
allow an exact estimation. Obviously, this condition is fulfilled
if matrix $C$ is equal to the identity.
\end{remarque}

Below, the synthesis of the observer matrices $L_x$  for the TEG of
figure \ref{figureGET} is given:
\begin{eqnarray*}
L_x=\begin{pmatrix} (4\gamma)^\ast & 6(4\gamma)^\ast  \\
\gamma^2(4\gamma)^\ast & 6\gamma^2(4\gamma)^\ast \\
\varepsilon & (3 \gamma)^\ast \end{pmatrix}
\end{eqnarray*}
Assumptions of proposition \ref{PropositionPente} being fulfilled,
it can easily be checked, by using toolbox Minmaxgd (see \cite{lhommeau00a}), that
$\sigma_\infty(x_i)=\sigma_\infty(\hat{x}_i)~\forall i \in [1,n]$
and that $Cx=C\hat{x} ~\forall (u,w)$ according to corollary
\ref{CorollaryOutputEgality}.

\section{Conclusion}
This paper\footnote{The authors are grateful to V. Reverdy for her valuable
linguistic help} has proposed a methodology  to design an observer for
($max,+$) linear systems.  The observer matrix is obtained thanks
to the residuation theory and is optimal in the sense that it is
the greatest which achieves the objective. It allows to compute a
state estimation lower than or equal to the real state and ensures
that the estimated output is equal to the system output. As a
perspective, this state estimation may be used in state feedback
control strategies as proposed in \cite{cottenceau01a,maia05}, and
an application to fault detection for manufacturing systems may be
envisaged. Furthermore, in order to deal with uncertain systems an extension can be
envisaged by considering interval analysis as it is done in
\cite{lhommeau04},\cite{hardouin09b} and more recently in \cite{diloreto09}.

\bibliographystyle{plain}        

\end{document}